\documentclass[9pt,technote]{IEEEtran}
\usepackage[cmex10]{amsmath}
\newtheorem{theorem}{Theorem}
\newtheorem{lemma}{Lemma}
\newtheorem{corollary}{Corollary}

\title{Reduction-Based Robustness Analysis of Linear Predictor Feedback for Distributed Input Delays}
\author{Anton~Ponomarev%
\thanks{Manuscript received August 2, 2014; revised January 5, 2015, May 5, 2015.}%
\thanks{A.~Ponomarev is with the Department of Control Theory, Saint Petersburg State University, Russia (e-mail: anton.pon.math@gmail.com).}%
\thanks{}}

\begin{document}

\maketitle

\begin{abstract}
Lyapunov--Krasovskii approach is applied to parameter- and delay-robustness analysis of the feedback suggested by Manitius and Olbrot for a linear time-invariant system with distributed input delay. A functional is designed based on Artstein's system reduction technique. It depends on the norms of the reduction-transformed plant state and original actuator state. The functional is used to prove that the feedback is stabilizing when there is a slight mismatch in the system matrices and delay values between the plant and controller.
\end{abstract}


\IEEEpeerreviewmaketitle

\IEEEpubid{0000--0000/00\$00.00~\copyright~2015 IEEE}

\section{Notation} We write $M>0$ or $M\geq 0$ to state that a symmetric real matrix $M$ is positive definite or positive semidefinite, respectively. Also in this case $\lambda_{\min}(M)$ and $\lambda_{\max}(M)$ represent the minimal and maximal eigenvalues of $M$. Vector norms being used are $\Vert x \Vert = \sqrt{x^T x}$ and $\Vert x\Vert_M = \sqrt{x^TMx}$, where $M>0$. Euclidean matrix norm is $\Vert M \Vert$.

The symbol $PC(T, X)$ stands for the space of piecewise continuous functions mapping $T\subset R$ into a Euclidean space $X$. The $L^2$ norm of $\varphi\in PC\big([-h,0), R^r\big)$ is $\Vert\varphi\Vert$, i.e.,
\begin{equation}
\Vert \varphi \Vert^2 =
	\int_{-h}^0 \Vert\varphi(\theta)\Vert^2 \,d\theta.
\end{equation}

Given $u\in PC(R, R^r)$, let $u_t$ be a function defined as $u_t(\theta)=u(t+\theta)$ for all $\theta\in[-h,0)$. The constant $h$ is specified below.

\section{Introduction}

\subsection{The problem} \label{SS_problem}

Consider the time-invariant system
\begin{equation} \label{plant0}
\dot{x}(t) = Ax(t) + \sum_{i=1}^{N} B_i u(t-h_i) + \int_{-h_\text{int}}^0 B_\text{int}(\theta) u(t+\theta) \,d\theta,
\end{equation}
where $x\in R^n$, $u\in R^r$, $h_1\geq 0$, $h_\text{int} \geq 0$, and $B_\text{int} \in PC\big([-h_\text{int}, 0], R^{n\times r}\big)$. For brevity, we will use Stieltjes integral notation and write the system under consideration as
\begin{equation} \label{plant}
\dot{x}(t) = Ax(t) + \int_{-h}^0 d\beta(\theta)u(t+\theta),
\end{equation}
where $h\geq \max\lbrace h_1, h_2, \dots, h_N, h_\text{int} \rbrace$,
\begin{equation}
\beta(\theta) = \sum_{i=1}^{N} B_i \chi(\theta+h_i)
	+ \int_{-h_\text{int}}^{\max\lbrace\theta,-h_\text{int}\rbrace} B_\text{int}(\tau) \,d\tau,
\end{equation}
and $\chi$ is the Heaviside step function.

The following control law was proposed for (\ref{plant}) in \cite{mani1979}:
\begin{equation} \label{control}
u(t) = Fx(t) + F\int_{-h}^0\int_\tau^0 e^{A(\tau-\theta)}\,d\beta(\tau) u(t+\theta) \,d\theta,
\end{equation}
where $F$ is a constant matrix. The feedback (\ref{control}) is called a \emph{predictor feedback} because it employs the plant's model (i.e., the matrices $A$ and $\beta$) to, in a sense, predict the future state of the plant.

Our goal is to investigate robustness of the feedback (\ref{control}). In terms of (\ref{plant0}), we are interested in:
\begin{itemize}
\item parametric robustness (small uncertainty in $A$, $B_i$, and $B_\text{int}$);
\item delay-robustness (small uncertainty in $h_i$ and $h_\text{int}$).
\end{itemize}

\subsection{Previous results overview} A range of methods is known to be suitable for analysis of linear systems of the form (\ref{plant}), (\ref{control}). Let us separate them into those using Lyapunov--Krasovskii functional analysis and those doing otherwise.

Most of the progress with non-Lyapunov techniques has been achieved in the area of systems with \emph{one discrete delay}, e.g., delay-robustness of a predictive controller \cite{mich2003,palm1980}, robustness with respect to a finite-sum implementation \cite{mond2002,melc2010}, robustness of an adaptive controller in presence of a disturbance \cite{eves2003}, and delay-robustness of a linear time-varying predictor feedback \cite{kara2013}. Furthermore, it has been shown in \cite{mond2003} that robustness with respect to a finite-sum implementation may be ensured by including a low-pass element in the control loop.

Lyapunov--Krasovskii analysis of systems with \emph{one discrete delay} was shown to succeed in proving delay-robustness of the predictor feedback \cite{krst2008} and robustness with respect to uncertain parameters \cite{jank2009}. Adaptive controllers were designed in \cite{bres2010,bres2012}. Recently, a predictor feedback for retarded \cite{khar2014a} and neutral \cite{khar2015} systems with state delays and an input delay was proposed, the closed loop's exponential stability being proven with a functional as well.

Lyapunov--Krasovskii analysis has been performed for \textit{distributed delays} too but less extensively. The results are closed-loop exponential stability \cite{beki2011} and stability with respect to an additive disturbance \cite{maze2012}. \textit{This paper} expands the list with parameter- and delay-robustness of the feedback.

\IEEEpubidadjcol

\subsection{Summary of the note} In Section III, the loop (\ref{plant}), (\ref{control}) is turned into (\ref{plant1}), (\ref{control1}) using the transformation (\ref{y}) borrowed from \cite{arts1982}. The Lyapunov--Krasovskii functional (\ref{func}) is then constructed. It includes the norms of the transformed plant state $y(t)$ and original actuator state $u_t$. Lemmas \ref{L_upbound}--\ref{L_lowbound} prove that this functional is quadratically bounded.

In Section IV, we show that a mismatch in $A$ and $\beta$ between (\ref{plant}) and (\ref{control}) introduces a distributed delay into the controller part of otherwise delay-free transformed system (\ref{plant1}), (\ref{control1}) which now becomes (\ref{plant1}), (\ref{control_rob1}). The delay does not significantly affect the system behavior if the mismatch is negligible. It leads to the \textit{main result}: closed-loop stability is robust (Theorem \ref{T_rob}). In order to facilitate a comparison of our approach with the preceding ones, we provide Corollary \ref{corol} together with its concise proof which is the special case of Theorem \ref{T_rob} for systems with one discrete delay.

\section{Lyapunov functional} Following Artstein \cite{arts1982}, let us introduce a new variable
\begin{equation} \label{y}
y(t) = x(t) + \int_{-h}^0 Q(\theta) u(t+\theta) \,d\theta,
\end{equation}
where
\begin{equation}
Q(\theta) = \int_{-h}^\theta e^{A(\tau-\theta)}\,d\beta(\tau).
\end{equation}
The closed loop (\ref{plant}), (\ref{control}) in the new variables takes the form
\begin{align}
\dot{y}(t) &= Ay(t) + Q(0)u(t), \label{plant1}\\
u(t) &= Fy(t). \label{control1}
\end{align}

It has been demonstrated by direct calculation in \cite{mani1979} that the eigenvalues of the closed loop (\ref{plant}), (\ref{control}) coincide with the eigenvalues of the matrix $A+Q(0)F$. Suppose the matrix is Hurwitz, so that the nominal closed loop is exponentially stable.

To apply Lyapunov approach, we aim at finding a functional $v(x,\varphi)$ defined for all $x\in R^n$ and $\varphi\in PC\big([-h,0),R^r\big)$ that admits upper and lower bounds proportional to $\Vert x\Vert^2 + \Vert \varphi\Vert^2$. We first construct a functional of this kind for the transformed loop (\ref{plant1}), (\ref{control1}) and then come back to the original variables.

Let us choose arbitrary matrices $W'>0$ and $W''>0$. Suppose then that $V>0$ is the solution of
\begin{equation} \label{lyapeq}
(A+Q(0)F)^T V + V (A+Q(0)F) = - (W' + 2F^T W'' F).
\end{equation}
For the loop (\ref{plant1}), (\ref{control1}), we propose the functional
\begin{equation}
\tilde{v}(y, \varphi) = \Vert y \Vert_V^2 +
	\int_{-h}^{\,0} e^{\sigma\theta} \Vert \varphi(\theta) \Vert_{W''}^2 \,d\theta
\end{equation}
defined for all $y\in R^n$ and $\varphi\in PC\big([-h,0),R^r\big)$, where
\begin{equation} \label{sigma}
\sigma = \frac{\lambda_{\min}(W')}{\lambda_{\max}(V)}.
\end{equation}
In the original variables it is
\begin{align} \label{func}
v(x, \varphi) = \left\Vert x + \!
	\int_{-h}^0 \! Q(\theta)\varphi(\theta) \,d\theta \right\Vert_V^2
	\! + \! \int_{-h}^{\,0} \!e^{\sigma\theta} \Vert \varphi(\theta) \Vert_{W''}^2 \,d\theta.
\end{align}

The following two lemmas prove that the functional (\ref{func}) has the required properties (upper and lower bounds).

\begin{lemma} \label{L_upbound}
\begin{equation}
v(x, \varphi) \leq M \big( \Vert x \Vert^2 + \Vert \varphi \Vert^2 \big),
\end{equation}
where
\begin{multline}
M = \max\big\lbrace
	2 \lambda_{\max}(V) \big\Vert e^{Ah} \big\Vert^2, \\
\lambda_{\max}(W'') +
	2h \lambda_{\max}(V) \max_{\theta\in[0,h]} \Vert Q(\theta) \Vert^2
\big\rbrace.
\end{multline}
\end{lemma}

\begin{IEEEproof}
One obtains this from (\ref{func}) using the Young's inequality, the triangle inequality, and the Cauchy--Schwarz inequality.
\end{IEEEproof}

\begin{lemma} \label{L_lowbound}
\begin{equation}
v(x, \varphi) \geq m_u \Vert \varphi \Vert^2, \quad
v(x, \varphi) \geq m_x \Vert x \Vert^2,
\end{equation}
where
\begin{gather}
m_u = e^{- \sigma h} \lambda_{\min}(W''), \\
m_x = \lambda_{\min}\Big( \big( V^{-1} + m_u^{-1} G \big)^{-1} \Big), \\
G = \int_{0}^{\,h} Q(\theta) Q^T(\theta) \,d\theta.
\end{gather}
\end{lemma}

Before we proceed to the proof, take notice that the lower bound of the functional in terms of the full norm follows from Lemma \ref{L_lowbound}:
\begin{equation}
v(x, \varphi) \geq \frac{\min\lbrace m_x, m_u \rbrace}{2}
	\big( \Vert x \Vert^2 + \Vert \varphi \Vert^2 \big).
\end{equation}
We leave the lemma as it is, though, lest the estimations be more conservative than necessary.

\begin{IEEEproof} The initial step is to estimate
\begin{equation} \label{lowbound1}
v(x, \varphi) \geq \Vert y \Vert_V^2 + m_u \Vert \varphi \Vert^2,
\end{equation}
where $y$ is linked to $x$ and $\varphi$ via (\ref{y}):
\begin{equation}
y = x + \int_{-h}^0 Q(\theta) \varphi(\theta) \,d\theta.
\end{equation}
The first inequality is obtained by dropping $\Vert y\Vert_V^2$ in (\ref{lowbound1}).

The idea used to establish the second inequality is to consider an optimal control problem: minimize the right hand side of (\ref{lowbound1}) with respect to $\varphi$.
Any function $\varphi\in PC\big( [-h,0), R^r \big)$ allows the decomposition (orthogonal projection on the rows of $Q$)
\begin{equation}
\varphi(\theta) = Q^T(\theta)c + \psi(\theta),
\end{equation}
where $c$ is a constant vector and $\psi\in PC\big( [-h,0), R^r \big)$ satisfies
\begin{equation}
\int_{-h}^0 Q(\theta)\psi(\theta)\,d\theta = 0.
\end{equation}
This representation lets one write
\begin{align}
y &= x + Gc,\\
\Vert \varphi \Vert^2 &= c^T G c + \Vert \psi \Vert^2 \geq c^T G c\\
\Rightarrow v(x, \varphi) &\geq \Vert Gc \Vert_V^2 + m_u c^T G c + 2 c^T GV x + \Vert x \Vert_V^2. \label{lowbound2}
\end{align}
The minimum of the quadratic estimation is at
\begin{equation}
c = - \big(G+m_uV^{-1}\big)^{-1}x.
\end{equation}
We substitute this $c$ into (\ref{lowbound2}) to find that
\begin{equation}
v(x,\varphi) \geq x^T \big( V^{-1} + m_u^{-1}G \big)^{-1} x,
\end{equation}
which leads to the desired inequality.
\end{IEEEproof}

\section{Robustness analysis}

\subsection{General result}

The controller's robustness with respect to a mismatch in the prediction model is analyzed here. In this scenario, the exact plant model (\ref{plant0}) or (\ref{plant}) is unknown but its estimation
\begin{equation} \label{plant0_rob}
\dot{x}(t) = \hat{A}x(t) + \sum_{i=1}^{N} \hat{B}_i u(t-\hat{h}_i) + \int_{-\hat{h}_\text{int}}^0 \hat{B}_\text{int}(\theta) u(t+\theta) \,d\theta
\end{equation}
or, equivalently,
\begin{equation} \label{plant_rob}
\dot{x}(t) = \hat{A}x(t) + \int_{-h}^0 d\hat{\beta}(\theta)u(t+\theta)
\end{equation}
is available. Here $\hat{h}_i\geq 0$, $\hat{h}_\text{int}\geq 0$, $B_\text{int} \in PC\big([-h_\text{int}, 0], R^{n\times r}\big)$, $h\geq \max\lbrace \hat{h}_1, \hat{h}_2, \dots, \hat{h}_N, \hat{h}_\text{int} \rbrace$,
\begin{equation}
\hat{\beta}(\theta) = \sum_{i=1}^{N} \hat{B}_i \chi(\theta+\hat{h}_i)
	+ \int_{-\hat{h}_\text{int}}^{\max\lbrace\theta,-\hat{h}_\text{int}\rbrace} \hat{B}_\text{int}(\tau) \,d\tau,
\end{equation}
and $\chi$ is the Heaviside step function. Observe that $h$ is the same in the nominal system (\ref{plant}) and its approximation (\ref{plant_rob}): for that, it is sufficient to take
\begin{equation}
h \geq \max\lbrace h_1, h_2, \dots, h_N, h_\text{int}, \hat{h}_1, \hat{h}_2, \dots, \hat{h}_N, \hat{h}_\text{int} \rbrace.
\end{equation}

The controller designed from (\ref{plant0_rob}) or (\ref{plant_rob}) would be
\begin{equation} \label{control_rob}
u(t) = F\left( x(t) + \int_{-h}^0 \hat{Q}(\theta) u(t+\theta) \,d\theta \right),
\end{equation}
where
\begin{equation}
\hat{Q}(\theta) = \int_{-h}^\theta
	e^{\hat{A}(\tau-\theta)} \,d\hat{\beta}(\tau).
\end{equation}
Regarding the choice of $F$, we demand that $A+Q(0)F$ be Hurwitz. However, matrices $A$ and $Q(0)$ are not known exactly due to parametric uncertainties. Nevertheless, suppose that one may establish some boundaries on $A$ and $Q(0)$ and choose a value of $F$ which renders $A+Q(0)F$ Hurwitz for all possible values of $A$ and $Q(0)$. This problem is not in the focus of the paper, so we take such $F$ as a given and assume hereafter that $A+Q(0)F$ is Hurwitz indeed.

After the transformation (\ref{y}), the plant is still (\ref{plant1}) and controller (\ref{control_rob}) is written as
\begin{equation} \label{control_rob1}
u(t) = F\left( y(t) +
	\int_{-h}^0 \Delta Q(\theta) u(t+\theta) \,d\theta
\right),
\end{equation}
where $\Delta Q(\theta) = \hat{Q}(\theta) - Q(\theta)$. One can clearly see how control delay reappears in the transformed loop (\ref{plant1}), (\ref{control_rob1}) due to imperfect modeling.

Let $v(t)$ be the value that the functional (\ref{func}) takes on a specific solution of the closed loop (\ref{plant}), (\ref{control_rob}).

\begin{lemma} \label{L_deriv_rob}
Along the solutions of the closed loop (\ref{plant}), (\ref{control_rob}) the functional (\ref{func}) satisfies
\begin{equation}
\dot{v}(t) \leq -\hat{\sigma} v(t),
\end{equation}
where
\begin{gather}
\hat{\sigma} = \sigma - k_1 \Vert \Delta Q \Vert - k_2 \Vert \Delta Q \Vert^2, \\
\Vert \Delta Q \Vert^2 = \int_{-h}^0 \Vert \Delta Q(\theta) \Vert^2 \,d\theta, \\
k_1 = \frac{\Vert VQ(0)F \Vert}{\min\lbrace\lambda_{\min}(V), m_u\rbrace}, \quad
k_2 = \frac{2\lambda_{\max}(W'') \Vert F\Vert^2}{m_u},
\end{gather}
$\sigma$ is (\ref{sigma}), and $m_u$ comes from Lemma \ref{L_lowbound}.
\end{lemma}

\begin{IEEEproof}
Differentiating $v(t)$, we use (\ref{lyapeq}) to get
\begin{align}
\dot{v}(t) &\leq 2y^T(t)V \big( Ay(t)+Q(0)u(t) \big) \nonumber \\
&\phantom{\leq}\:- \sigma \int_{-h}^{\,0} e^{\sigma\theta} \Vert u(t+\theta) \Vert_{W''}^2 \,d\theta
	+ \Vert u(t)\Vert_{W''}^2 \\
&\leq -\sigma v(t) + 2 y^T(t)VQ(0)F \int_{-h}^0 \Delta Q(\theta) u(t+\theta) \,d\theta \nonumber \\
&\phantom{\leq}\:+ 2 \Vert u(t) - Fy(t)\Vert_{W''}^2.
\end{align}
Further estimations include
\begin{align}
2 y^T(t)VQ(0) & F\int_{-h}^0 \Delta Q(\theta) u(t+\theta) \,d\theta \nonumber \\
&\leq 2 \Vert VQ(0)F \Vert \Vert \Delta Q \Vert
	\Vert y(t) \Vert \Vert u_t \Vert \\
&\leq \Vert VQ(0)F \Vert \Vert \Delta Q \Vert
	\big( \Vert y(t) \Vert^2 + \Vert u_t \Vert^2 \big) \\
&\leq k_1 \Vert \Delta Q \Vert \,v(t)
\end{align}
and
\begin{align}
2\Vert u(t) - Fy(t)\Vert_{W''}^2 &\leq 2\lambda_{\max}(W'') \Vert F\Vert^2
	\big(\Vert\Delta Q\Vert \Vert u_t\Vert\big)^2\\
&\leq k_2 \Vert\Delta Q\Vert^2 \,v(t).
\end{align}
In the end, we arrive at the estimation $\dot{v}(t) \leq -\hat{\sigma} v(t)$.
\end{IEEEproof}

The following is our main result.

\begin{theorem} \label{T_rob}
If $A+Q(0)F$ is Hurwitz and (\ref{plant0_rob}) approximates (\ref{plant0}) closely enough that \textit{(see the remark after the proof)}
\begin{equation} \label{rob}
\Vert \Delta Q \Vert^2 < \left(
	\frac{\sqrt{k_1^2 + 4k_2\sigma} - k_1}{2k_2}
\right)^2,
\end{equation} 
then the closed loop (\ref{plant}), (\ref{control_rob}) is exponentially stable:
\begin{align}
\Vert x(t) \Vert^2 &\leq
	\frac{M}{m_x} \, e^{-\hat{\sigma} t} \big( \Vert x(0) \Vert^2 + \Vert u_0 \Vert^2 \big), \\
\Vert u_t \Vert^2 &\leq
	\frac{M}{m_u} \, e^{-\hat{\sigma} t} \big( \Vert x(0) \Vert^2 + \Vert u_0 \Vert^2 \big).
\end{align}
Here $\Vert\Delta Q\Vert$, $k_1$, $k_2$, $M$, $m_x$, $m_u$, and $\hat{\sigma}$ are defined in Lemmas \ref{L_upbound}--\ref{L_deriv_rob}, and we mean that (\ref{plant0_rob}) approximates (\ref{plant0}) in the sense that
\begin{align}
\hat{A} &\approx A, \\
\hat{B}_i &\approx B_i \text{ for all } i=1,2,\dots,N, \\
\hat{B}_\text{int} &\approx B_\text{int} \text{ uniformly on }
	\big[-\min\lbrace h_\text{int}, \hat{h}_\text{int} \rbrace, 0 \big], \\
\hat{h}_i &\approx h_i \text{ for all } i=1,2,\dots,N, \\
\hat{h}_\text{int} &\approx h_\text{int}.
\end{align}
\end{theorem}

\begin{IEEEproof}
It follows from Lemmas \ref{L_upbound}--\ref{L_deriv_rob}.
\end{IEEEproof}

\textit{Remark:} Let us explain why $\Vert\Delta Q\Vert$ becomes small when approximation of (\ref{plant0}) with (\ref{plant0_rob}) is almost perfect. If all parameters of (\ref{plant0_rob}) closely match those of (\ref{plant0}), then $\hat{\beta}$ and $\beta$ by definition are uniformly close on $[-h, 0]$ except for a finite number of small intervals like $[-h_i, -\hat{h}_i)$, $[-\hat{h}_i, -h_i)$, $[-h_\text{int}, -\hat{h}_\text{int})$ or $[-\hat{h}_\text{int}, -h_\text{int})$. It follows then that $\Delta Q$ is uniformly small on $[-h, 0]$ except the small intervals where it is bounded (and the bound is independent of the intervals' lengths). Smallness of $\Vert\Delta Q\Vert$ ensues.

See the proof of Corollary \ref{corol} for a quantitative example of how $\Vert \Delta Q \Vert$ depends on a mismatch in the value of a discrete delay.

\subsection{Special case: one discrete delay}

To highlight the conditions imposed by Theorem \ref{T_rob} on discrete delay values, we supplement it with an application to the widely studied single discrete delay case.

\begin{corollary} \label{corol}
The closed loop
\begin{align}
\dot{x}(t) &= Ax(t) + Bu(t-\delta), \label{plantrobex} \\
u(t) &= F\left( x(t) +
	\int_{-\hat{\delta}}^0 e^{-A(\hat{\delta}+\theta)}
	B u(t+\theta) \,d\theta \right), \label{controlrobex}
\end{align}
where $\hat{\delta}\geq 0$, $\delta\geq 0$, is exponentially stable if there exist $V>0$, $W'>0$, $W''>0$ such that
\begin{equation}
\big(A+e^{-A\delta}BF\big)^T V + V \big(A+e^{-A\delta}BF\big) = - \big(W' + 2F^T W'' F\big)
\end{equation}
and
\begin{equation}
|\hat{\delta}-\delta| < \frac{1}{\Vert A\Vert} \ln \left(
	1 + \frac{\sqrt{1 + k_3 e^{2\Vert A\Vert h}}-1}
	{e^{2\Vert A\Vert h}} \right),
\end{equation}
where
\begin{gather}
k_3 = \frac{2\Vert A\Vert}{\Vert B\Vert^2}
	\left( \frac{\sqrt{k_1^2 + 4k_2\sigma} - k_1}{2k_2} \right)^2, \\
k_1 = \frac{\Vert V e^{-A\delta}B F \Vert}{\min\left\lbrace \lambda_{\min}(V),
	e^{-\sigma h}\lambda_{\min}(W'') \right\rbrace}, \\
k_2 = \frac{2\lambda_{\max}(W'') \Vert F\Vert^2}{e^{-\sigma h} \lambda_{\min}(W'')}, \\
\sigma = \frac{\lambda_{\min}(W')}{\lambda_{\max}(V)}, \quad
h = \max\lbrace \delta, \hat{\delta} \rbrace.
\end{gather}
If $A=0$, then the delay mismatch condition is
\begin{equation}
|\hat{\delta}-\delta| < \left( \frac{\sqrt{k_1^2+4k_2\sigma} - k_1}{2\Vert B\Vert k_2} \right)^2.
\end{equation}
\end{corollary}

\begin{IEEEproof}
Given a solution $\big( x(t), u_t \big)$ of the closed loop (\ref{plantrobex}), (\ref{controlrobex}), we study the behavior of
\begin{equation}
v(t) = \Vert y(t) \Vert_V^2 +
	\int_{-h}^0 e^{\sigma\theta}
	\Vert u(t+\theta) \Vert_{W''}^2 \,d\theta,
\end{equation}
where
\begin{gather}
y(t) = x(t) + \int_{-h}^0 Q(\theta) u(t+\theta) \,d\theta, \label{yrob} \\
Q(\theta) = \begin{cases}
	e^{-A(\delta+\theta)} B, &\text{if } \theta\in[-\delta, 0],\\
	0, &\text{if } \theta<-\delta.
\end{cases}
\end{gather}
The control law (\ref{controlrobex}) is equivalently represented as
\begin{equation}
u(t) = F\left(
	y(t) + \int_{-h}^0 \Delta Q(\theta) u(t+\theta) \,d\theta
\right),
\end{equation}
where $\Delta Q(\theta) = \hat{Q}(\theta) - Q(\theta)$ and
\begin{equation}
\hat{Q}(\theta) = \begin{cases}
	e^{-A(\hat{\delta}+\theta)} B, &\text{if } \theta\in[-\hat{\delta}, 0],\\
	0, &\text{if } \theta<-\hat{\delta}.
\end{cases}
\end{equation}

Differentiating $v(t)$ and performing some estimations explained in the proof of Lemma \ref{L_deriv_rob}, we get
\begin{equation}
\dot{v}(t) \leq \big(-\sigma + k_1 \Vert \Delta Q \Vert +
	k_2 \Vert \Delta Q \Vert^2 \big) \,v(t),
\end{equation}
where
\begin{equation}
\Vert \Delta Q \Vert^2 = \int_{-h}^0 \Vert \Delta Q(\theta) \Vert^2 \,d\theta.
\end{equation}
If
\begin{equation} \label{raterob}
-\sigma + k_1 \Vert \Delta Q \Vert + k_2 \Vert \Delta Q \Vert^2 < 0,
\end{equation}
then exponential stability of (\ref{plantrobex}), (\ref{controlrobex}) is guaranteed by Lemmas \ref{L_upbound}--\ref{L_deriv_rob}.

Let us estimate $\Vert \Delta Q \Vert$. If $h=\hat{\delta}>\delta$,
\begin{equation}
\Delta Q(\theta) = \begin{cases}
	e^{-A(h+\theta)} B,
		&\text{if } \theta\in[-h, -\delta), \\
	\left(I - e^{A|\hat{\delta}-\delta|}\right)
		e^{-A(\hat{\delta}+\theta)} B,
		&\text{if } \theta\in[-\delta, 0].
\end{cases}
\end{equation}
$\Vert \Delta Q \Vert$ is small because $\Delta Q(\theta)$ is bounded on the first interval which has small length, and small on the second one which length is bounded:
\begin{align}
\Vert \Delta Q \Vert^2 & = \int_{-h}^{-\delta}
	\Vert \Delta Q(\theta) \Vert^2 \,d\theta +
	\int_{-\delta}^0 \Vert \Delta Q(\theta) \Vert^2 \,d\theta \\
& \leq \frac{\Vert B\Vert^2}{2\Vert A\Vert}
	\left( e^{2\Vert A\Vert |\hat{\delta}-\delta|} - 1 \right) \nonumber \\
& \phantom{=}\: + \frac{\Vert B\Vert^2}{2\Vert A\Vert}
	\left( e^{2\Vert A\Vert h} - 1 \right)
	\left( e^{\Vert A\Vert |\hat{\delta}-\delta|} - 1 \right)^2. \label{DQ}
\end{align}
If $A = 0$, then by continuity (\ref{DQ}) turns into
\begin{equation} \label{DQ_spec}
\Vert \Delta Q \Vert^2 \leq \Vert B\Vert^2 |\hat{\delta}-\delta|.
\end{equation}
If $h=\delta>\hat{\delta}$,
\begin{equation}
\Delta Q(\theta) = \begin{cases}
	-e^{-A(h+\theta)} B,
		&\text{if } \theta\in[-h, -\hat{\delta}), \\
	\left(I - e^{-A|\hat{\delta}-\delta|} \right)
		e^{-A(\hat{\delta}+\theta)} B,
		&\text{if } \theta\in[-\hat{\delta}, 0],
\end{cases}
\end{equation}
but the same estimation (\ref{DQ}) or (\ref{DQ_spec}) holds true. The corollary's premise then resolves in (\ref{raterob}), thus proving exponential stability.
\end{IEEEproof}

It should be mentioned that Corollary \ref{corol} is not the best result available for the single discrete delay case as our estimations are quite conservative. The problem has been widely studied with more accurate results achieved, e.g., in \cite{mich2003,kara2013}.

\section{Conclusions} A Lyapunov--Krasovskii functional has been devised for a linear time-invariant system with distributed input delay closed by a predictor feedback. The functional is based on Artstein's system reduction technique.

The framework is shown to be convenient for robustness analysis of the closed loop: it alleviates the proof of closed-loop exponential stability when controller's predictive model (matrices and delay values) differs slightly from the actual plant.

Future research may include a generalization of the approach for nonlinear systems and systems with both state and control delays.

\section*{Acknowledgment} The author thanks Prof.~A.~Zhabko, Prof.~M.~Krstic, and anonymous reviewers for valuable feedback and advice.

\bibliographystyle{IEEEtran}
\bibliography{IEEEabrv,paper_lyapfun_bib}

\end{document}